\documentclass[article,11pt]{amsart}

\usepackage{amsfonts}
\usepackage{amsmath}
\usepackage{mathrsfs}
\usepackage{graphicx}
\usepackage{color}
\usepackage{epsfig}
\usepackage[font={small}]{caption}
\usepackage{xcolor}

\numberwithin{equation}{section}
\theoremstyle{plain}
\newtheorem{thm}{Theorem}[section]
\newtheorem{rem}{Remark}[section]

\newtheorem{lem}{Lemma}[section]

\newcommand{\dE}{\mathbb{E}}
\newcommand{\dR}{\mathbb{R}}
\newcommand{\dS}{\mathbb{S}}

\newcommand{\dP}{\mathbb{P}}

\newcommand{\dQ}{\mathbb{Q}}
\newcommand{\dZ}{\mathbb{Z}}

\newcommand{\cN}{\mathcal{N}}

\newcommand{\rI}{\mathrm{I}}
\newcommand{\cF}{\mathcal{F}}

\newcommand{\cM}{\mathcal{M}}

\newcommand{\veps}{\varepsilon}

\newcommand{\ind}{\mbox{1}\kern-.25em \mbox{I}}
\font\calcal=cmsy10 scaled\magstep1

\def\build#1_#2^#3{\mathrel{\mathop{\kern 0pt#1}\limits_{#2}^{#3}}}
\def\liml{\build{\longrightarrow}_{}^{{\mbox{\calcal L}}}}

\def\videbox{\mathbin{\vbox{\hrule\hbox{\vrule height1.4ex \kern.6em\vrule height1.4ex}\hrule}}}
\def\demend{\hfill $\videbox$\\}

\keywords{Elephant random walk; Center of mass; Multi-dimensional martingales; Almost sure convergence; Asymptotic normality}

\begin{document}
\title[On the center of mass of the elephant random walk]
{On the center of mass of the elephant random walk \vspace{1ex}}
\author{Bernard Bercu}
\thanks{The corresponding author is Bernard Bercu, email address: bernard.bercu@math.u-bordeaux.fr}
\dedicatory{\normalsize University of Bordeaux, France}
\address{Universit\'e de Bordeaux, Institut de Math\'ematiques de Bordeaux,
UMR 5251, 351 Cours de la Lib\'eration, 33405 Talence cedex, France.}
\email{bernard.bercu@math.u-bordeaux.fr}
\email{lucile.laulin@math.u-bordeaux.fr}
\author{Lucile Laulin}


\begin{abstract}
Our goal is to investigate the asymptotic behavior of the center of mass of the elephant random walk, which is a discrete-time random walk on integers
with a complete memory of its whole history. In the diffusive and critical regimes, we establish the almost sure convergence, the law of iterated logarithm 
and the quadratic strong law for the center of mass of the elephant random walk. The asymptotic normality of the center of mass, properly normalized, 
is also provided. Finally, we prove a strong limit theorem for the center of mass in the superdiffusive regime. All our analysis relies on asymptotic results 
for multi-dimensional martingales.
\end{abstract}
\maketitle


\section{Introduction}
\label{S-I}
Let $(S_n)$ be a standard random walk in $\dR^d$. The center of mass $G_n$ of $S_n$ is defined by
\begin{equation}
\label{CM}
G_n=\frac{1}{n}\sum_{k=1}^n S_k.
\end{equation}
The question of the asymptotic behavior of $G_n$ was first raised by Paul Erd\"{o}s. Very recently, 
Lo and Wade \cite{LoWade18} extended the results of Grill \cite{Grill88}
by studying the asymptotic behavior of $(G_n)$. More precisely, let $S_n=X_1+\cdots+X_n$ where the increments $(X_n)$ 
are independent and identically distributed square integrable random vectors of $\dR^d$ with mean $\mu$ and covariance matrix
$\Gamma$. They proved the strong law of large numbers
\begin{equation}
\label{CMSL}
\lim_{n \rightarrow \infty} \frac{1}{n}G_n=\frac{1}{2}\mu \hspace{1cm} \text{a.s.}
\end{equation}
together with the asymptotic normality,
\begin{equation}
\label{CMAN}
\frac{1}{\sqrt{n}} \Bigl( G_n - \frac{n}{2}\mu \Bigr) \liml \cN \Bigl(0, \frac{1}{3}\Gamma \Bigr).
\end{equation}
Curiously, no other references are availabe on the asymptotic behavior of the center of mass.
The proofs of many results on the convex hull $C_n$ as well as on the center of mass $G_n$
rely on independence and exchangeability of the increments of the walk. For example, one can observe that
\begin{equation}
\label{CMEXP}
G_n= \frac{1}{n}\sum_{k=1}^n S_k=\frac{1}{n} \sum_{k=1}^n \bigl( n-k+1 \bigr)X_k
\end{equation}
shares the same distribution as
$$
\Sigma_n=\frac{1}{n} \sum_{k=1}^n k X_k.
$$
A natural question concerns the asymptotic behavior of $G_n$ in other situations where the increments of the walk are not
independent and not identically distributed. Moreover, geometrical features of the random walk paths of $(S_n)$ is a
subject of ongoing interest. For example, the convex hull $C_n=\text{Conv}(S_1, \ldots, S_n)$ of the $n$ first steps
of $(S_n)$ have recently received renewed attention \cite{Kabluchko17},  \cite{Vysotsky18}. 
More particularly, for the random walk in $\dR^2$, the strong law of large numbers and the asymptotic normality of
the perimeter and the diameter of $C_n$ were established in \cite{McRedmond18}, \cite{Wade15}.

\ \vspace{1ex}\\
In this paper, we investigate the asymptotic behavior of the center of mass of the multi-dimensional elephant random walk.
It is a fascinating discrete-time random walk on $\dZ^d$ where $d \geq 1$, which has a complete memory of its whole history.
The increments depend on all the past of the walk and they are not exchangeable.
The elephant random walk (ERW) was introduced by Sch\"utz and Trimper \cite{Schutz04} in the early 2000s, in order to investigate how
long-range memory affects the random walk and induces a crossover from a diffusive
to superdiffusive behavior. It was referred to as the ERW in allusion to the traditional saying that elephants can always remember where they have been before.
The elephant starts at the origin at time zero, $S_0 = 0$. At time $n = 1$,  the elephant moves in one of the $2d$ directions with the same probability $1/2d$.  
Afterwards, at time $n+1$, the elephant chooses uniformly at random an integer $k$ among the previous times $1,\ldots,n$. Then, it moves exactly in 
the same direction as that of time $k$ with probability $p$ or in one of the $2d-1$ remaining directions with the same probability $(1-p)/(2d-1)$, where the parameter $p$ 
stands for the memory parameter of the ERW \cite{Bercu19}. Therefore, the position of the elephant at time $n+1$ is given by
\begin{equation} 
\label{POSMERW}
S_{n+1} = S_n + X_{n+1}
\end{equation}
where $X_{n+1}$ is the $(n+1)$-th increment of the random walk.
The ERW shows three differents regimes depending on the location of its memory parameter $p$
with respect to the critical value
\begin{equation}
\label{DEFPD}
p_d=\frac{2d+1}{4d}.
\end{equation}
A wide literature is now available on the ERW in dimesion $d=1$ where $p_d=3/4$. A strong law of large numbers
and a central limit theorem for the position $S_n$, properly normalized, were established 
in the diffusive regime $p< 3/4$ and the critical regime $p=3/4$, see \cite{Baur16}, \cite{Coletti17}, \cite{ColettiN17}, \cite{Schutz04}
and the recent contributions \cite{BercuHG19}, \cite{Bertoin2020}, \cite{Businger18}, \cite{Coletti19}, \cite{Fan2020}, \cite{Takei2020}, \cite{Vazquez19}. 
The superdiffusive regime $p>3/4$ is much harder to handle.
Bercu \cite{Bercu18} proved that the limit of the position of the ERW is not Gaussian. Quite recently, Kubota and Takei \cite{Kubota19} showed that 
the fluctuation of the ERW around its limit in the superdiffusive regime is Gaussian. Finally, Bercu and Laulin in \cite{Bercu19}
extended all the results of \cite{Bercu18} to the multi-dimensional ERW where $d \geq 1$. 
\vspace{1ex} \\
Our strategy for proving asymptotic results for the center of mass of the elephant random walk (CMERW) is as follows.
On the one hand, the behavior of position $S_n$ is closely related to the one of the 
sequence $(M_n)$ defined, for all $n \geq 0$, by $M_n = a_nS_n$ with
$a_1 = 1$ and, for all $n \geq 2$,
\begin{equation}
\label{DEFAN}
a_n=\prod_{k=1}^{n-1}\Bigl(\frac{k}{k+a} \Bigr)= \frac{\Gamma(a+1)\Gamma(n)}{\Gamma(n+a)}
\end{equation}
where $\Gamma$ stands for the Euler Gamma function and $a$ is the fundamental parameter of the ERW
defined by
\begin{equation}
\label{DEFA}
a=\frac{2dp-1}{2d-1}.
\end{equation}
It follows from the very definition of the ERW that at any time $n \geq 1$, 
\begin{equation}
\dE\left[X_{n+1}|\cF_n \right]  =  \Bigl( p I_d - \frac{1-p}{2d-1} I_d \Bigr) \frac{1}{n} \sum_{k=1}^nX_{k} =
\frac{1}{n}\Bigl(\frac{2dp-1}{2d-1}\Bigr) S_n = \frac{a}{n} S_n
\hspace{1cm}\text{a.s.}
\label{CEX}
\end{equation}
Hence, we obtain from \eqref{POSMERW} and \eqref{CEX} that for any $n \geq 1$,
\begin{equation}
\label{CES1}
\dE\left[S_{n+1}|\cF_n\right] = \Bigl(  1+\frac{a}{n} \Bigr) S_n 
\hspace{1cm}\text{a.s.}
\end{equation}
Therefore, we deduce from \eqref{DEFAN} and \eqref{CES1} that
\begin{equation}
\label{CEM1}
\dE\left[M_{n+1}|\cF_n\right] = a_{n+1}\Bigl(  1+\frac{a}{n} \Bigr) S_n=a_nS_n=M_n 
\hspace{1cm}\text{a.s.}
\end{equation}
It means that $(M_n)$ is a locally square-integrable martingale adapted to the filtration $(\cF_n)$
where $\cF_n=\sigma(X_1, \ldots,X_n)$. It can be rewritten  \cite{Bercu19} in the additive form 
\begin{equation}
\label{DEFMN}
M_n=\sum_{k=1}^n a_k \veps_k 
\end{equation}
where $\veps_1=S_1$ and, for all $n \geq 2$,
\begin{equation}
\label{DEFVEPSN}
\veps_n = S_n - \Bigl(\frac{a_{n-1}}{a_n}\Bigr)S_{n-1}= S_n - \Bigl(1+\frac{a}{n-1}\Bigr)S_{n-1}.
\end{equation}
On the other hand, an analogue of equation \eqref{CMEXP} is given by
\begin{eqnarray}
G_n  &=& \frac{1}{n} \sum_{k=1}^n S_k =   \frac{1}{n} \sum_{k=1}^n  \frac{1}{a_k}M_k =  \frac{1}{n} \sum_{k=1}^n \frac{1}{a_k}\sum_{\ell=1}^k a_\ell\veps_\ell 
=  \frac{1}{n} \sum_{k=1}^n a_k \veps_k 
\sum_{\ell=k}^n \frac{1}{a_\ell}, \nonumber \\
&=& \frac{1}{n}\sum_{k=1}^n a_k(b_n-b_{k-1})\veps_k
\label{DEFGNADD}
\end{eqnarray}
where the sequence $(b_n)$ is given by $b_0=0$ and, for all $n \geq 1$,
\begin{equation}
\label{DEFBN}
b_n=\sum_{k=1}^n \frac{1}{a_k}.
\end{equation} 
Denoting
\begin{equation}
\label{DEFNN}
N_n=\sum_{k=1}^n a_k b_{k-1}\veps_k,
\end{equation}
it is straightforward to see that 
$\dE\left[N_{n+1}|\cF_n\right] =N_n$ a.s.
since $\dE\left[\veps_{n+1}|\cF_n\right] =0$. Hence, 
$(N_n)$ is also a locally square-integrable martingale adapted to the filtration $(\cF_n)$.
We deduce from \eqref{DEFGNADD} that
\begin{equation} 
\label{EMCM}
  G_n= \frac{1}{n}(b_n M_n-N_n).
\end{equation}
Relation \eqref{EMCM} allows us to establish the asymptotic behavior of the CMERW via
an extensive use of the strong law of large numbers and the central limit theorem for multi-dimensional martingales
\cite{Chaabane00}, \cite{Duflo97}, \cite{Hall80}, \cite{Touati91}.
\ \vspace{4ex} \\
The paper is organized as follows. The main results are given in Section \ref{S-MR}. 
We first investigate the diffusive regime $p<p_d$ and we establish the almost sure convergence, the law of iterated logarithm and the quadratic strong law for the CMERW. The asymptotic normality of the CMERW, properly normalized, is also provided. Next, we prove similar results in the critical regime $p=p_d$. Finally, we 
establish a strong limit theorem in the superdiffusive regime $p>p_d$. 
our martinagle approach is described in Section \ref{S-MA} while
all technical proofs are postponed to Appendices A, B and C.


\section{Main results}
\label{S-MR}


\subsection{The diffusive regime}
Our first result deals with the strong law of large numbers for the CMERW in the diffusive regime where $0 \leq p < p_d$.
The following strong law for the CMERW will be deduced as a simple consequence of the strong law for $(S_n)$.

\begin{thm}
\label{T-ASCV-CM-DR}
We have the almost sure convergence
\begin{equation}
\label{ASCV-CM-DR1}
 \lim_{n \rightarrow \infty} \frac{1}{n}G_n=0 \hspace{1cm} \text{a.s.}
\end{equation}
\end{thm}
\begin{rem}
For any $\alpha>1/2$, we have the more precise result 
$$
\lim_{n \to \infty} \frac{1}{n^{\alpha}}G_n = 0 \hspace{1cm}\text{a.s.}
$$
\end{rem}

\begin{figure}[ht]
\vspace{-2ex}
\begin{center}
\includegraphics[scale=0.3]{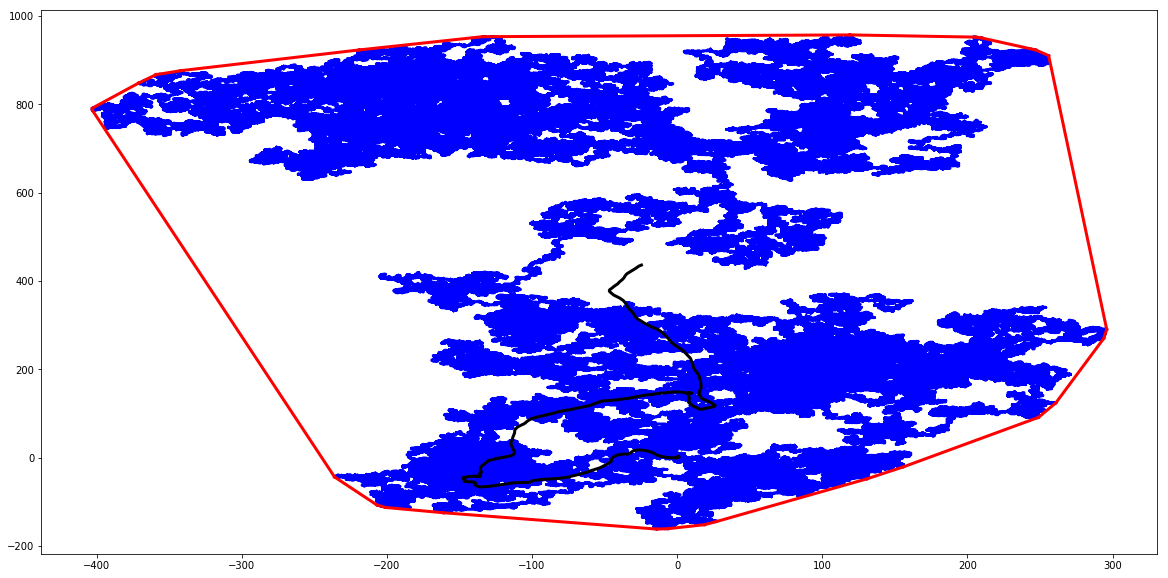}
\vspace{-2ex}
\caption{The 2-dimensional ERW in blue, the CMERW in black and the convex hull in red, for $n=10^6$ steps and 
a diffusive memory parameter $p=1/2$.}
\label{Fig-DR}
\end{center}
\end{figure}

\noindent
The almost sure rates of convergence for CMERW are as follows.
\begin{thm}
\label{T-LFQLIL-CM-DR}
We have the quadratic strong law
\begin{equation}
\label{LFQ-CM-DR1}
 \lim_{n \rightarrow \infty} \frac{1}{\log n} \sum_{k=1}^n \frac{1}{k^2}G_k G_k^T=\frac{2}{3(1-2a)(2-a)d}I_d \hspace{1cm} \text{a.s.}
\end{equation}
where $I_d$ stands for the identity matrix of order $d$.
In particular,
\begin{equation}
\lim_{n \rightarrow \infty} \frac{1}{ \log n} \sum_{k=1}^n  \frac{\|G_k\|^2}{k^2}=  \frac{2}{3(1-2a)(2-a)} \hspace{1cm} \text{a.s.}
 \label{LFQ-CM-DR2}
\end{equation}
Moreover, we also have the upper-bound in the law of iterated logarithm
\begin{equation}
 \limsup_{n \rightarrow \infty} \frac{\|G_n\|^2}{2 n \log \log n} \leq  \frac{\bigl(\sqrt{3}+\sqrt{1-2a}\bigr)^2}{3(a+1)^2(1-2a)d} \hspace{1cm} \text{a.s.}  
\label{UPLIL-CM-DR3}
\end{equation}
\end{thm}
\noindent
We are now interested in the asymptotic normality of the CMERW.

\begin{thm}
\label{T-CLT-CM-DR}
We have the asymptotic normality
\begin{equation}
\label{CLT-CM-DR1}
\frac{1}{\sqrt{n}}G_n \liml \cN \Bigl(0, \frac{2}{3(1-2a)(2-a)d}I_d \Bigr).
\end{equation}
\end{thm}

\begin{rem}
One can observe from Theorem 3.3 in \cite{Bercu19} that the ratio of the asymptotic variances between the CMERW 
and the ERW is given by
\begin{equation*}
	R(a)=\frac{2}{3(2-a)}.
\end{equation*} 
In the diffusive regime, this ratio lies between $2/9$ and $4/9$ and it is always smaller than $1$, as one can see in Figure \ref{Fig-DR}.
Moreover, in the special case where the elephant moves in one of the $2d$ directions with the same probability $p=1/2d<p_d$, 
it follows from \eqref{DEFA} that the fundamental parameter $a=0$. Consequently, we deduce from
\eqref{CLT-CM-DR1} that
\begin{equation*}
\frac{1}{\sqrt{n}}G_n \liml \cN \Bigl(0,\frac{1}{3d}I_d\Bigr).
\end{equation*}
We find again the asymptotic normality \eqref{CMAN} 
where the mean value $\mu=0$ and the covariance matrix $\Gamma=\frac{1}{d}I_d$.
\end{rem}


\subsection{The critical regime}
Hereafter, we investigate the critical regime where the memory parameter $p=p_d$.

\begin{thm}
\label{T-ASCV-CM-CR}
We have the almost sure convergence
\begin{equation}
\label{ASCV-CM-CR1}
 \lim_{n \rightarrow \infty} \frac{1}{\sqrt{n} \log n}G_n=0 \hspace{1cm} \text{a.s.}
\end{equation}
\end{thm}
\begin{rem}
For any $\alpha>1/2$, we have the more precise result 
$$
\lim_{n \to \infty} \frac{1}{\sqrt{n}(\log n)^{\alpha}}G_n = 0 \hspace{1cm}\text{a.s.}
$$
\end{rem}

\begin{figure}[ht]
\vspace{-1ex}
\begin{center}
\includegraphics[scale=0.3]{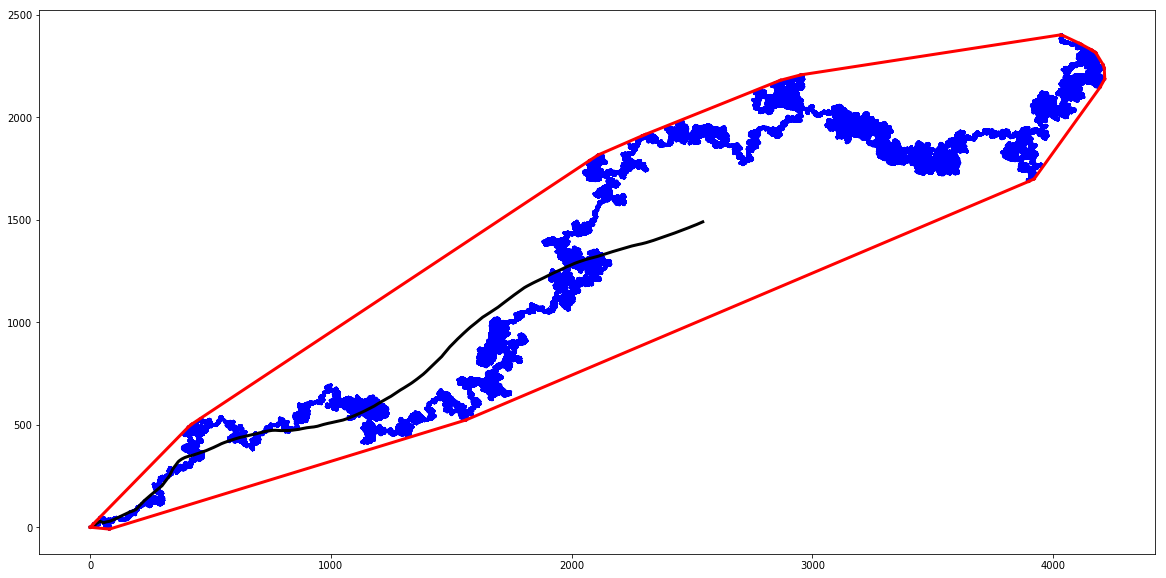}
\vspace{-2ex}
\caption{The 2-dimensional  ERW in blue, the CMERW in black and the convex hull in red, for $n=10^6$ steps and 
a critical memory parameter $p=5/8$.}
\label{Fig-CR}
\end{center}
\end{figure}

\noindent
The almost sure rates of convergence for the CMERW are as follows.
\begin{thm}
\label{T-LFQLIL-CM-CR}
We have the quadratic strong law
\begin{equation}
\label{LFQ-CM-CR1}
 \lim_{n \rightarrow \infty} \frac{1}{\log \log n} \sum_{k=2}^n  \frac{1}{(k \log k)^2}G_k G_k^T=  \frac{4}{9d}I_d \hspace{1cm} \text{a.s.}
\end{equation}
In particular,
\begin{equation}
\lim_{n \rightarrow \infty} \frac{1}{\log \log n} \sum_{k=2}^n  \frac{\|G_k\|^2}{(k \log k)^2}=  \frac{4}{9} \hspace{1cm} \text{a.s.}
 \label{LFQ-CM-CR2}
\end{equation}
Moreover, we also have the law of iterated logarithm
\begin{equation}
 \limsup_{n \rightarrow \infty} \frac{\|G_n\|^2}{2 n \log n \log \log \log n}=  \frac{4}{9d} \hspace{1cm} \text{a.s.}
 \label{LIL-CM-CR3}
\end{equation}
\end{thm}
\noindent
Our next result concerns the asymptotic normality of the CMERW.

\begin{thm}
\label{T-CLT-CM-CR}
We have the asymptotic normality
\begin{equation}
\label{CVL-CM-CR1}
\frac{1}{\sqrt{n\log n}}G_n \liml \cN \left(0, \frac{4}{9d}I_d \right).
\end{equation}
\end{thm}

\begin{rem}
In the critical regime, the ratio of the asymptotic variances between the CMERW and the ERW is $4/9$.
\end{rem}


\subsection{The superdiffusive regime}
Finally, we focus our attention on the superdiffusive regime where 
$p>p_d$. The almost sure convergence of $(S_n)$, properly normalized by $n^a$, yields the following
strong limit theorem for the CMERW.

\begin{thm}
\label{T-ASCV-CM-SR}
We have the almost sure convergence
\begin{equation}
\label{ASCV-CM-SR1}
 \lim_{n \rightarrow \infty} \frac{1}{n^a}G_n=G \hspace{1cm} \text{a.s.}
\end{equation}
where the limiting value G is a non-degenerate random vector of $\dR^d$.
Moreover, we also have the mean square convergence
\begin{equation}
\label{ASCV-CM-SR2}
\lim_{n \to \infty} \dE\Bigl[\Bigl\|\frac{1}{n^a}G_n-G\Bigr\|^2\Bigr]=0.
\end{equation}
\end{thm}

\begin{rem}
The expected value of  $G$ is zero and its covariance matrix is given by
\begin{equation*}
\label{COVG}
\dE\left[GG^T\right]=\frac{1}{d(a+1)^2(2a -1)^2\Gamma(2a-1)} I_d.
\end{equation*}
The distribution of $G$ is sub-Gaussian but far from being known.
\end{rem}

\begin{figure}[ht]
\vspace{-1ex}
\begin{center}
\includegraphics[scale=0.3]{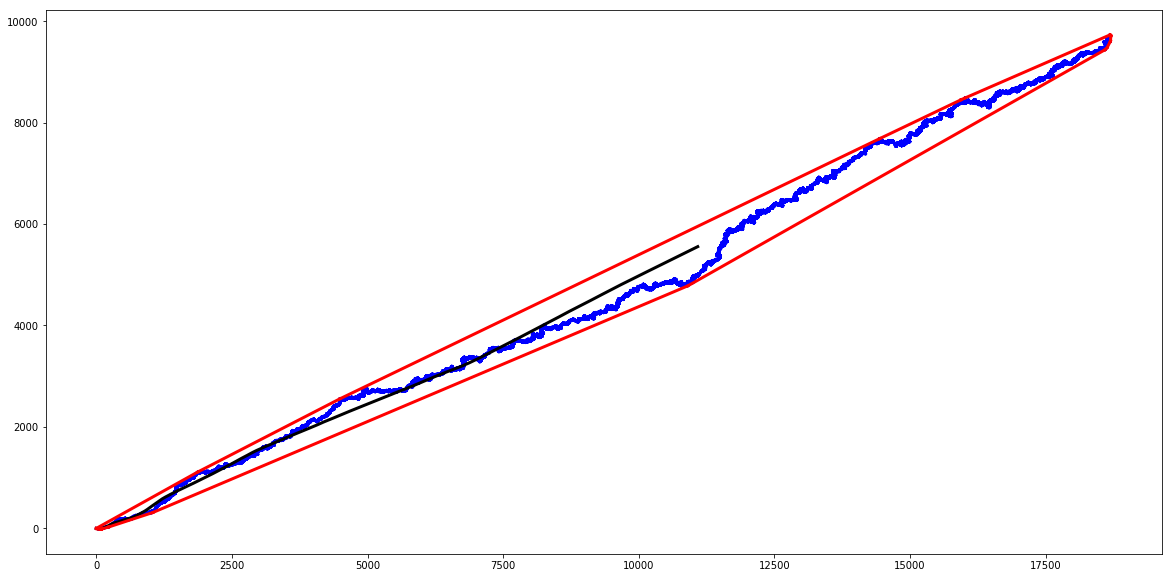}
\vspace{-2ex}
\caption{The 2-dimensional  ERW in blue, the CMERW in black and the convex hull in red, for $n=10^6$ steps and 
a superdiffusive memory parameter $p=3/4$.}
\label{Fig-SR}
\end{center}
\end{figure}


\section{A multi-dimensional martingale approach}
\label{S-MA}


We already saw from \eqref{EMCM} that the CMERW can be rewritten as
\begin{equation*} 
G_n= \frac{1}{n}(b_n M_n-N_n).
\end{equation*}
In  order to investigate the asymptotic behavior of $(G_n)$, we introduce
the multi-dimensional martingale $(\cM_n)$ defined by
\begin{equation}
\label{DEFMNN}
 \cM_n = \begin{pmatrix}M_n \\ N_n\end{pmatrix}
\end{equation} 
where $(M_n)$ and $(N_n)$ are the two locally square-integrable martingales given by \eqref{DEFMN} and \eqref{DEFNN}. 
The main difficulty we face here is that the predictable quadratic variation of $(M_n)$ and $(N_n)$ increase
to infinity with two different speeds. A matrix normalization is necessary to establish 
the asymptotic behavior of the CMERW. Let $(V_n)$ be the sequence of positive definite diagonal matrices of order $2d$ given by
\begin{equation}
\label{DEFVN}
V_n= \frac{1}{n\sqrt{n}}
\begin{pmatrix}
b_n & 0 \\ 0 & 1 
\end{pmatrix}\otimes I_d
\end{equation}
where $A\otimes B$ stands for the Kronecker product of the matrices $A$ and $B$.

\begin{lem}
\label{L-CVGIP}
The sequence $(\cM_n)$ is a locally square-integrable martingale of $\dR^{2d}$.  
Its predictable quadratic variation $\langle \cM \rangle_n$ satisfies in the diffusive regime where $a<1/2$, 
\begin{equation}
\label{CVGIPMN}
\lim_{n \rightarrow \infty} V_n \langle \cM \rangle_n V_n^T = V \hspace{1cm} \text{a.s.}
\end{equation}
where the limiting matrix
\begin{equation}
\label{DEFMATRIXV}
V= \frac{1}{d(a+1)^2}\begin{pmatrix} 
                  \frac{1}{1-2a} & \frac{1}{2-a} \\
                  \frac{1}{2-a} & \frac{1}{3}   
                  \end{pmatrix}\otimes I_d.
\end{equation}
\end{lem}

\begin{rem}
Via the same lines as in the proof of Lemma \ref{L-CVGIP}, we find that in the critical regime $a=1/2$,
the sequence of normalization matrices $(V_n)$ has to be replaced by
\begin{equation}
\label{DEFWN}
W_n= \frac{1}{n\sqrt{n \log n}}
\begin{pmatrix}
b_n & 0 \\ 0 & 1 
\end{pmatrix}\otimes I_d.
\end{equation}
Moreover, the limiting matrix in \eqref{CVGIPMN} must be changed by
\begin{equation}
\label{DEFMATRIXW}
W= \frac{4}{9d}\begin{pmatrix} 
                  1 & 0 \\
                  0 & 0   
                  \end{pmatrix}\otimes I_d.
\end{equation}
\end{rem}

\begin{proof}
The increments of the ERW are bounded by $1$, that is for any time $n\geq 1$, $\|X_{n} \| = 1$. Hence,
it follows from \eqref{POSMERW} that $\|S_n\| \leq n$ and $\|G_n\| \leq n$ which imply that
$\| M_n \| \leq n a_n$ and $\| N_n \| \leq n a_n b_n + n^2$. 
We  already saw in Section \ref{S-I} that $(\cM_n)$ is a locally square-integrable martingale.
Denote $\Delta M_n=M_n-M_{n-1}$, and similarly for other processes. It follows from \eqref{DEFMN}, \eqref{DEFVEPSN} and \eqref{DEFNN}
that the predictable quadratic variation associated with
$(\cM_n)$ is the square matrix of order $2d$ given, for all $n \geq 1$, by
\begin{equation}
\label{IPMN1}
\langle \cM \rangle_n \!=\! \sum_{k=1}^{n} 
\dE \! \left[\!
\begin{pmatrix} \Delta M_k \\  \Delta N_k \end{pmatrix}
\!\!\begin{pmatrix} \Delta M_k \\  \Delta N_k \end{pmatrix}^{\!T}\!
\Bigl|  \cF_{k-1}\!\right]\!
\!=\! \sum_{k=1}^{n} a_k^2
\begin{pmatrix}1 & b_{k-1} \\ b_{k-1} &b_{k-1}^2 \end{pmatrix}
\!\otimes \dE \bigl[\veps_k \veps_k^T \bigl|  \cF_{k-1}\bigr].
\end{equation}
Moreover, we deduce from formulas (A.7) and (B.3) in \cite{Bercu19} that for all $n \geq 1$,
\begin{equation}
\dE\left[\veps_{n+1}\veps_{n+1}^T|\cF_n\right] = \frac{1}{d}I_d
+a\Bigl( \frac{1}{n}\Sigma_n - \frac{1}{d}I_d\Bigr) - \left( \frac{a}{n}\right) ^2 S_nS_n^T\hspace{1cm}  \text{a.s.}
\label{CEMEPS2}
\end{equation}
where $\Sigma_n$ is a random positive definite matrix of order $d$ satisfying
\begin{equation}
\label{CVGSIGMAN}
\lim_{n\to\infty}\frac{1}{n}\Sigma_n = \frac{1}{d}I_d \hspace{1cm} \text {a.s.}
\end{equation}
Consequently, we obtain from \eqref{IPMN1} together with \eqref{CEMEPS2} that
\begin{equation}
\langle \cM \rangle_n =  \frac{1}{d}\sum_{k=1}^n a_k^2\begin{pmatrix} 1 & b_{k-1}\\ b_{k-1} & b_{k-1}^2  \end{pmatrix} \otimes I_d 
        + a\sum_{k=1}^{n-1} a_{k+1}^2 \begin{pmatrix} 1 & b_{k}\\ b_{k} & b_{k}^2 \end{pmatrix}\otimes\Bigl( \frac{1}{k}\Sigma_k 
        - \frac{1}{d}I_d\Bigr) - \xi_n
\label{IPMN2}
\end{equation} 
where 
$$
\xi_n=a^2\sum_{k=1}^{n-1} \left(\frac{a_{k+1}}{k}\right)^2 \begin{pmatrix} 1 & b_{k}\\ b_{k} & b_{k}^2 \end{pmatrix} \otimes S_k S_k^T.
$$ 
According to Theorem 3.1 in \cite{Bercu19}, the remainder $\xi_n$ plays a negligible role as
\begin{equation}
\label{ASCVGSN}
\lim_{n \rightarrow \infty} \frac{S_n}{n}=0 \hspace{1cm} \text{a.s.}
\end{equation}
Hereafter, it is not hard to see that
$$
V_n \Bigl( \sum_{k=1}^n a_k^2\begin{pmatrix} 1 & b_{k-1}\\ b_{k-1} & b_{k-1}^2  \end{pmatrix} \otimes I_d \Bigr)V_n^T=
\frac{1}{n^3} \begin{pmatrix} b_n^2   \sum_{k=1}^n a_k^2 & b_n  \sum_{k=1}^n a_k^2b_{k-1}\\ b_n  \sum_{k=1}^n a_k^2b_{k-1} &  \sum_{k=1}^n a_k^2b_{k-1}^2  \end{pmatrix}
\otimes I_d.
$$
Furthermore, from a well-known property of the Euler Gamma function, we have
\begin{equation}
\label{CVGAMMAN}
\lim_{n \rightarrow \infty} \frac{\Gamma(n+a)}{\Gamma(n) n^a}= 1.
\end{equation}
Hence, we obtain from \eqref{DEFAN}, \eqref{DEFBN} and \eqref{CVGAMMAN} that
\begin{equation}
\label{CVGANBN}
\lim_{n \rightarrow \infty} n^a a_n= \Gamma(a+1) \hspace{1cm}\text{and}\hspace{1cm}
\lim_{n \rightarrow \infty} \frac{b_n}{n^{a+1}}= \frac{1}{\Gamma(a+2)}.
\end{equation}
Consequently, as soon as $a<1/2$, we immediately find from \eqref{CVGANBN} that
\begin{eqnarray*}
\lim_{n \rightarrow \infty}  \frac{b_n^2}{n^3} \sum_{k=1}^n a_k^2 & = & \frac{1}{(1-2a)(a+1)^2}, \\
\lim_{n \rightarrow \infty}  \frac{b_n}{n^3} \sum_{k=1}^n a_k^2 b_{k-1} & = & \frac{1}{(2-a)(a+1)^2}, \\
\lim_{n \rightarrow \infty}  \frac{1}{n^3} \sum_{k=1}^n  a_k^2 b_{k-1}^2 & = & \frac{1}{3(a+1)^2}.
\end{eqnarray*}
Therefore,
\begin{equation}
\label{IPMN3}
\lim_{n \rightarrow \infty} V_n \Bigl( \sum_{k=1}^n a_k^2\begin{pmatrix} 1 & b_{k-1}\\ b_{k-1} & b_{k-1}^2  \end{pmatrix} \otimes I_d \Bigr)V_n^T=
\frac{1}{(a+1)^2}\begin{pmatrix} 
                  \frac{1}{1-2a} & \frac{1}{2-a} \\
                  \frac{1}{2-a} & \frac{1}{3}   
                  \end{pmatrix}\otimes I_d.
\end{equation}
Finally, it follows from the combinaison of \eqref{CVGSIGMAN}, \eqref{IPMN2}, \eqref{ASCVGSN}  and \eqref{IPMN3} that
\begin{equation}
\label{IPMN4}
\lim_{n \rightarrow \infty} V_n \langle \cM \rangle_n V_n^T = \frac{1}{d(a+1)^2}\begin{pmatrix} 
                  \frac{1}{1-2a} & \frac{1}{2-a} \\
                  \frac{1}{2-a} & \frac{1}{3}   
                  \end{pmatrix}\otimes I_d. \hspace{1cm} \text{a.s.}
\end{equation}
which is exactly what we wanted to prove.
\end{proof}

\section*{Appendix A. Two non-standard results on martingales}
\renewcommand{\thesection}{\Alph{section}}
\renewcommand{\theequation}{\thesection.\arabic{equation}}
\setcounter{section}{1}
\setcounter{equation}{0}
\setcounter{thm}{0}
The proofs of our main results rely on two non-standard central limit theorem and quadratic strong law for multi-dimensional martingales.
A simplified version of Theorem 1 of Touati \cite{Touati91} is as follows.
\begin{thm}
\label{T-CLT}
Let $(\cM_n)$ be a locally square-integrable martingale of $\dR^\delta$ adapted to a filtration $(\cF_n)$,  
with predictable quadratic variation $\langle \cM \rangle_n$.
Let $(V_n)$ be a sequence of non-random square matrices of order $\delta$ such that $\| V_n \|$ 
decreases to $0$ as $n$ goes to infinity. Assume that there exists a symmetric and positive semi-definite matrix $V$ 
such that
\begin{equation*}
V_n \langle \cM \rangle_n V_n^T \overset{\displaystyle \dP}{\underset{n\to\infty}{\longrightarrow}} V. \leqno (\textnormal{H.1})
\end{equation*}
Moreover, assume that Lindeberg's condition is satisfied, that is for all $\veps >0$,
\begin{equation*}
 \sum_{k=1}^n \dE\bigl[\|V_n \Delta \cM_k \|^2 \rI_{\{\|V_n\Delta \cM_k \|>\veps\}}\bigl|\cF_{k-1}\bigr] \overset{\displaystyle \dP}{\underset{n\to\infty}{\longrightarrow}}   0 \leqno (\textnormal{H.2})
\end{equation*}
where $\Delta \cM_n=\cM_n - \cM_{n-1}$. Then, we have the asymptotic normality
\begin{equation}
\label{CLTM}
  V_n \cM_n \liml \cN \bigl(0, V \bigr).
\end{equation}
\end{thm}

\noindent
The quadratic strong law requires more restrictive assumptions. The following result is a simplified version of Theorem 2.1 of Chaabane and Maaouia \cite{Chaabane00}
where the normalization matrices $(V_n)$ are diagonal.
\begin{thm}
\label{T-QSL}
Let $(\cM_n)$ be a locally square-integrable martingale of $\dR^\delta$ adapted to a filtration $(\cF_n)$,  
with predictable quadratic variation $\langle \cM \rangle_n$.
Let $(V_n)$ be a sequence of non-random positive definite diagonal matrices of order $\delta$ such that
its diagonal terms decrease to zero at polynomial rates. Assume that $(\textnormal{H.1})$ and $(\textnormal{H.2})$ hold almost surely.
Moreover, suppose that there exists $\beta\in]1,2]$ such that
\begin{equation*}
\sum_{n=1}^{\infty} \frac{1}{\bigl(\log  (\det V_{n}^{-1})^2\bigr)^{\beta}}\dE\bigl[\|V_{n} \Delta \cM_{n}\|^{2\beta}\bigl|\cF_{n-1}\bigr]<\infty \hspace{1cm} \text{a.s.}
\leqno (\textnormal{H.3})
\end{equation*}
Then, we have the quadratic strong law
\begin{equation}
\label{LFQM}
\lim_{n \rightarrow \infty}  \frac{1}{\log(\det V_n^{-1})^2}\sum_{k=1}^n
\Bigl(\frac{(\det V_{k})^2 - (\det V_{k+1})^2}{(\det V_k)^2}\Bigr)V_k\cM_k\cM_k^T V_k^T = V \hspace{1cm} \text{a.s.}
\end{equation}
\end{thm}

\section*{Appendix B. Proofs of the almost sure convergence results}
\renewcommand{\thesection}{\Alph{section}}
\renewcommand{\theequation}{\thesection.\arabic{equation}}
\setcounter{section}{2}
\setcounter{equation}{0}
\setcounter{thm}{0}

\subsection{The diffusive regime.}
\ \vspace{2ex}\\
\noindent{\bf Proof of Theorem \ref{T-ASCV-CM-DR}.}
We already saw from Theorem 3.1 in \cite{Bercu19} that 
\begin{equation}
\label{ASCVG-SN-DR}
\lim_{n \rightarrow \infty} \frac{S_n}{n}=0 \hspace{1cm} \text{a.s.}
\end{equation}
Consequently, the almost sure convergence \eqref{ASCV-CM-DR1} immediately follows from \eqref{ASCVG-SN-DR} together with the
Toeplitz lemma given e.g. by Lemma 2.2.13 in \cite{Duflo97}.
\demend
\ \vspace{1ex}\\
\noindent{\bf Proof of Theorem \ref{T-LFQLIL-CM-DR}.}
Our goal is to check that all the hypotheses of Theorem \ref{T-QSL} are satisfied. Thanks to Lemma \ref{L-CVGIP}, hypothesis
$(\textnormal{H.1})$ holds almost surely. In order to verify that 
Lindeberg's condition $(\textnormal{H.2})$ is satisfied, we have from \eqref{DEFMNN} together with \eqref{DEFMN}, \eqref{DEFNN} and $V_n$ given by \eqref{DEFVN}
that for all $1 \leq k \leq n$
$$
V_n \Delta \cM_k = \frac{a_k}{n \sqrt{n}} \begin{pmatrix} b_n \veps_k \\  b_{k-1} \veps_k\end{pmatrix},
$$
which implies that 
\begin{equation}
\label{VNMN-BOUND}
	\| V_n \Delta \cM_k \|^2 
 \leq \frac{2 a_k^2 b_n^2}{n^3}  \| \veps_k \|^2.
\end{equation}
Consequently, we obtain that for all $\veps >0$,
\begin{eqnarray}
 \sum_{k=1}^n \dE\bigl[\|V_n \Delta \cM_k \|^2 \rI_{\{\|V_n\Delta \cM_k \|>\veps\}}\bigl|\cF_{k-1}\bigr]
            & \leq & \frac{1}{\veps^2}\sum_{k=1}^n \dE\bigl[\|V_n \Delta \cM_k \|^4\bigl|\cF_{k-1}\bigr], \nonumber\\
            & \leq & \frac{4b_n^4}{\veps^2 n^6}\sum_{k=1}^n a_k^4\dE\bigl[\|\veps_k\|^4\bigl|\cF_{k-1}\bigr], \nonumber \\
            & \leq & \frac{4b_n^4}{\veps^2 n^6 } \sup_{1\leq k \leq n} \dE\bigl[\|\veps_k\|^4\bigl|\cF_{k-1}\bigr]\sum_{k=1}^n a_k^4.
\label{LINDEBERG-DR1} 
 \end{eqnarray}
However, it follows from the right-hand side of formula (4.11) in \cite{Bercu19} that 
\begin{equation}
\label{MAJMOMEPS}
\sup_{1\leq k \leq n} \dE\bigl[\|\veps_k\|^4\bigl|\cF_{k-1}\bigr] \leq \frac{4}{3} \hspace{1cm} \text{a.s.}
\end{equation}	
Therefore, we infer from \eqref{LINDEBERG-DR1} that for all $\veps >0$,
\begin{equation}
 \sum_{k=1}^n \dE\bigl[\|V_n \Delta \cM_k \|^2 \rI_{\{\|V_n\Delta \cM_k \|>\veps\}}\bigl|\cF_{k-1}\bigr]
             \leq  \frac{16 \,b_n^4}{3 \,\veps^2 n^6 } \sum_{k=1}^n a_k^4 \hspace{1cm} \text{a.s.}
\label{LINDEBERG-DR2} 
\end{equation}
Moreover, we have from \eqref{CVGANBN} that
\begin{equation}
\label{LINDEBERG-DR3}
  b_n^4 \sum_{k=1}^n a_k^4 = O( n^5 ).
\end{equation}
Consequently, \eqref{LINDEBERG-DR2} together with \eqref{LINDEBERG-DR3} ensure that
Lindeberg's condition $(\textnormal{H.2})$ holds almost surely, that is
for all $\veps >0$,
\begin{equation}
\label{LINDEBERG-DR4}
\lim_{n \rightarrow \infty} \sum_{k=1}^n \dE\bigl[\|V_n \Delta \cM_k \|^2 \rI_{\{\|V_n\Delta \cM_k \|>\veps\}}\bigl|\cF_{k-1}\bigr]= 0 \hspace{1cm} \text{a.s.}
\end{equation}
We will now check that condition $(\textnormal{H.3})$ is satisfied in the special case $\beta=2$, that is
\begin{equation}
\label{QSL-DR1}
\sum_{n=1}^{\infty} \frac{1}{\bigl(\log  (\det V_{n}^{-1})^2\bigr)^{2}}\dE\bigl[\|V_{n} \Delta \cM_{n}\|^{4}\bigl|\cF_{n-1}\bigr]<\infty \hspace{1cm} \text{a.s.}
\end{equation}
We have from \eqref{DEFVN} that
\begin{equation}
\label{DETVN}
      \det V_n^{-1} = \Bigl(\frac{n^{3/2}}{b_n}\Bigr)^d.
\end{equation} 
Hence, we find from \eqref{CVGANBN} and \eqref{DETVN} that
\begin{equation}
\label{LOGDETVN}
 \lim_{n \rightarrow \infty}  \frac{\log  (\det V_{n}^{-1})^2 }{\log n} = d(1-2a).
\end{equation} 
Consequently, we can replace $\log  (\det V_{n}^{-1})^2$ by $\log n$ in \eqref{QSL-DR1}. Hereafter, we obtain from \eqref{VNMN-BOUND} and \eqref{MAJMOMEPS} that
\begin{eqnarray}
\sum_{n=2}^{\infty} \frac{1}{(\log n)^2}\dE\bigl[\|V_{n} \Delta \cM_{n}\|^{4}\bigl|\cF_{n-1}\bigr]
       & = & O \Bigl( \sum_{n=1}^{\infty} \frac{1}{(\log n)^2}\frac{a_n^4 b_n^4}{n^6} \dE\bigl[\|\veps_n\|^4\bigl|\cF_{n-1}\bigr] \Bigr), \nonumber \\
       & = & O \Bigl( \sum_{n=1}^{\infty} \frac{1}{(\log n)^2}\frac{a_n^4 b_n^4}{n^6}  \Bigr). 
\label{QSL-DR2}
\end{eqnarray}
However, we have from \eqref{CVGANBN} that
\begin{equation}
\label{QSL-DR3}
\lim_{n \rightarrow \infty}  \frac{a_n^4b_n^4}{n^4} =  \frac{1}{(a+1)^4}.
\end{equation}
Therefore, \eqref{QSL-DR2} together with \eqref{QSL-DR3} immediately lead to \eqref{QSL-DR1}.
We are now in a position to apply the quadratic strong law given by Theorem \ref{T-QSL}. We have from
\eqref{LFQM} and \eqref{LOGDETVN} that
\begin{equation}
\label{QSL-DR4}
\lim_{n \rightarrow \infty}  \frac{1}{\log n}\sum_{k=1}^n
\Bigl(\frac{(\det V_{k})^2 - (\det V_{k+1})^2}{(\det V_k)^2}\Bigr)V_k\cM_k\cM_k^T V_k^T = d(1-2a) V \hspace{1cm} \text{a.s.}
\end{equation}
where the limiting matrix $V$ is given by \eqref{DEFMATRIXV}.
However, it follows from \eqref{DEFGNADD}, \eqref{DEFMNN} and \eqref{DEFVN} that
\begin{equation}
\label{DECOMPMARTGN}
\frac{1}{\sqrt{n}}G_n= v^T V_n \cM_n
\hspace{1cm}\text{where}\hspace{1cm}
v=\begin{pmatrix}1 \\ -1\end{pmatrix}\otimes I_d.
\end{equation} 
Consequently, we deduce from \eqref{QSL-DR4} and \eqref{DECOMPMARTGN} that
\begin{equation}
\label{QSL-DR5}
\lim_{n \rightarrow \infty}  \frac{1}{\log n}\sum_{k=1}^n
\Bigl(\frac{(\det V_{k})^2 - (\det V_{k+1})^2}{(\det V_k)^2}\Bigr)\frac{1}{k}G_kG_k^T = d(1-2a) v^TVv \hspace{1cm} \text{a.s.}
\end{equation}
Furthermore, we obtain from \eqref{DETVN} and  \eqref{CVGANBN} that
$$
\lim_{n \rightarrow \infty} n\Bigl(\frac{(\det V_{n})^2 - (\det V_{n+1})^2}{(\det V_n)^2}\Bigr)=d(1-2a).
$$
Hence, \eqref{QSL-DR5} clearly leads to convergence \eqref{LFQ-CM-DR1},
\begin{equation}
\label{QSL-DR6}
\lim_{n \rightarrow \infty}  \frac{1}{\log n}\sum_{k=1}^n
\frac{1}{k^2}G_kG_k^T =  v^TVv = \frac{2}{3(1-2a)(2-a)d}I_d \hspace{1cm} \text{a.s.}
\end{equation}
By taking the trace on both sides of \eqref{QSL-DR6}, we also obtain that
\begin{equation}
 \lim_{n \rightarrow \infty}  \frac{1}{\log n}\sum_{k=1}^n \frac{\|G_k\|^2}{k^2} =  \frac{2}{3(1-2a)(2-a)}\hspace{1cm}\text{a.s.}
\end{equation}
Finally, we shall proceed to the proof of the upper-bound \eqref{UPLIL-CM-DR3} in the law of iterated logarithm. 
Denote
 \begin{equation}
\label{DEFTAUN}
  \tau_n = \sum_{k=1}^n a_k^2 b_{k-1}^2.
\end{equation}
We already saw from \eqref{QSL-DR3} that $a_{n}^4b_{n-1}^4 \tau_n^{-2}$ is equivalent to $9n^{-2}$. It implies
that
\begin{equation}
\label{CONDLIL-DR}
\sum_{n=1}^{+\infty} \frac{a_{n}^4b_{n-1}^4}{\tau_n^2} < \infty.
\end{equation}
Moreover, we have from \eqref{CVGSIGMAN}, \eqref{IPMN2}, \eqref{ASCVGSN} and \eqref{DEFTAUN} that 
$$
 \lim_{n \rightarrow \infty}  \frac{1}{\tau_n} \langle N \rangle_n=\frac{1}{d} I_d \hspace{1cm}\text{a.s.}
$$
Consequently, we deduce from the law of iterated logarithm for martingales due to Stout \cite{Stout74}, see also Corollary 6.4.25 in \cite{Duflo97}, 
that $(N_n)$ satisfies for any vector $u \in \dR^d$,
\begin{eqnarray}
 \limsup_{n \rightarrow \infty} \Bigl(\frac{1}{2 \tau_n \log \log \tau_n}\Bigr)^{1/2} \langle u, N_n \rangle & = &
 -\liminf_{n \rightarrow \infty} \Bigl(\frac{1}{2 \tau_n \log \log \tau_n}\Bigr)^{1/2} \langle u, N_n \rangle \nonumber \\
 & = & \frac{1}{\sqrt{d}}\|u\| \hspace{1cm} \text{a.s.}
 \label{LIL-CM-DR1}
\end{eqnarray}
However, since $\tau_n$ is equivalent to $n^3/3(a+1)^2$, \eqref{LIL-CM-DR1} immediately lead to
\begin{eqnarray}
 \limsup_{n \rightarrow \infty} \Bigl(\frac{1}{2 n \log \log n}\Bigr)^{1/2} \frac{1}{n}\langle u, N_n \rangle & = &
 -\liminf_{n \rightarrow \infty} \Bigl(\frac{1}{2 n \log \log n}\Bigr)^{1/2} \frac{1}{n}\langle u, N_n \rangle \nonumber \\
 & = & \frac{1}{\sqrt{3d}(a+1)}\|u\| \hspace{1cm} \text{a.s.}
 \label{LIL-CM-DR2}
\end{eqnarray}
Furthermore, it was already shown by formula (5.17) in \cite{Bercu19} that for any vector $u \in \dR^d$,
\begin{eqnarray}
 \limsup_{n \rightarrow \infty} \Bigl(\frac{n^{2a}}{2 n \log \log n}\Bigr)^{1/2} \langle u, M_n\rangle & = & 
 -\liminf_{n \rightarrow \infty} \Bigl(\frac{n^{2a}}{2 n \log \log n}\Bigr)^{1/2} \langle u, M_n\rangle \nonumber \\
 & = & \frac{\Gamma(a+1)}{\sqrt{d(1-2a)}}\|u\| \hspace{1cm} \text{a.s.}
\label{LIL-CM-DR3}
\end{eqnarray}
Therefore, we deduce from \eqref{EMCM} and \eqref{CVGANBN} together with \eqref{LIL-CM-DR2} 
and \eqref{LIL-CM-DR3} that for any vector $u$ of $\dR^d$,
\begin{equation*}
\limsup_{n \rightarrow \infty} \Bigl(\frac{1}{2 n \log \log  n}\Bigr)^{1/2} \langle u, G_n \rangle = 
\limsup_{n \rightarrow \infty} \Bigl(\frac{1}{2 n \log \log  n}\Bigr)^{1/2} \frac{1}{n}\langle u, b_n M_n-N_n \rangle
\vspace{-2ex}
\end{equation*}
\begin{eqnarray}
\hspace{1cm}  & \leq &
\limsup_{n \rightarrow \infty} \Bigl(\frac{1}{2 n \log \log  n}\Bigr)^{1/2} \frac{1}{n}\langle u, b_n M_n\rangle  + 
\limsup_{n \rightarrow \infty} \Bigl(\frac{1}{2 n \log \log  n}\Bigr)^{1/2} \frac{1}{n}\langle u, -N_n \rangle\nonumber \\
& \leq &
\limsup_{n \rightarrow \infty} \Bigl(\frac{1}{2 n \log \log  n}\Bigr)^{1/2} \frac{1}{n}\langle u, b_n M_n \rangle  - 
\liminf_{n \rightarrow \infty} \Bigl(\frac{1}{2 n \log \log  n}\Bigr)^{1/2} \frac{1}{n}\langle u, N_n \rangle\nonumber \\
& \leq &
 \frac{\|u\|}{\sqrt{d}(a+1)}\Bigl(\frac{1}{\sqrt{1-2a}} +\frac{1}{\sqrt{3}}\Bigr)
 \hspace{1cm} \text{a.s.}
 \label{LIL-CM-DR4}
\end{eqnarray}
By the same token, we also find that for any vector $u$ of $\dR^d$,
\begin{equation*}
\liminf_{n \rightarrow \infty} \Bigl(\frac{1}{2 n \log \log  n}\Bigr)^{1/2} \langle u, G_n \rangle = 
\liminf_{n \rightarrow \infty} \Bigl(\frac{1}{2 n \log \log  n}\Bigr)^{1/2} \frac{1}{n}\langle u, b_n M_n-N_n \rangle
\vspace{-2ex}
\end{equation*}
\begin{eqnarray}
\hspace{1cm}  & \geq &
 - \frac{\|u\|}{\sqrt{d}(a+1)}\Bigl(\frac{1}{\sqrt{1-2a}} +\frac{1}{\sqrt{3}}\Bigr)
 \hspace{1cm} \text{a.s.}
 \label{LIL-CM-DR5}
\end{eqnarray}
Consequently, we obtain from \eqref{LIL-CM-DR4}  and \eqref{LIL-CM-DR5} that for any vector $u$ of $\dR^d$,
\begin{equation}
 \label{LIL-CM-DR6}
\limsup_{n \rightarrow \infty} \Bigl(\frac{1}{2 n \log \log  n}\Bigr) \langle u, G_n \rangle^2 \leq  
\frac{\|u\|^2}{d(a+1)^2}\Bigl(\frac{1}{\sqrt{1-2a}} +\frac{1}{\sqrt{3}}\Bigr)^2
 \hspace{1cm} \text{a.s.}
 \end{equation}
One can observe that the upper-bound in \eqref{LIL-CM-DR6} is close to the optimal bound
$$
(vu)^TVvu= \frac{2\|u\|^2}{3(1-2a)(2-a)d}.
$$
Finally, by taking all rational points on the unit sphere $\dS^{d-1}$ in $\dR^d$, the bound in \eqref{LIL-CM-DR6} holds simultaneously for all of them,
which implies that
$$
 \limsup_{n \rightarrow \infty} \frac{\| G_n\|^2}{2 n  \log \log n} \leq   \!\!
 \sup_{u \in \dQ^d \cap \dS^{d-1}} \!\! \limsup_{n \rightarrow \infty} \frac{\langle u, G_n \rangle^2}{2 n  \log \log n}  \leq 
\frac{\bigl(\sqrt{3}+\sqrt{1-2a}\bigr)^2}{3(a+1)^2(1-2a)d} 
 \hspace{1cm} \text{a.s.}
$$
completing the proof of Theorem \ref{T-LFQLIL-CM-DR}.
\demend

\subsection{The critical regime.}
\ \vspace{2ex}\\
\noindent{\bf Proof of Theorem \ref{T-ASCV-CM-CR}.}
We have from Theorem 3.4 in \cite{Bercu19} that 
\begin{equation}
\label{ASCVG-SN-CR}
\lim_{n \to \infty} \frac{1}{\sqrt{n}\log n}S_n = 0 \hspace{1cm} \text{a.s.}
\end{equation}
Hence,  \eqref{ASCV-CM-CR1} clearly follows from \eqref{ASCVG-SN-CR} together with the
Toeplitz lemma.
\demend

\noindent{\bf Proof of Theorem \ref{T-LFQLIL-CM-CR}.}
The proof of the quadratic strong law \eqref{LFQ-CM-CR1} is left to the reader as it 
follows essentially the same lines as that of \eqref{LFQ-CM-DR1}. The only minor change is that the
matrix $V_n$ has to be replaced by the matrix $W_n$ defined in \eqref{DEFMATRIXW}.
We shall now proceed to the proof of the law of iterated logarithm
given by \eqref{LIL-CM-CR3}. On the one hand, it follows from \eqref{LIL-CM-DR2} with $a=1/2$ that
for any vector $u \in \dR^d$,
\begin{eqnarray*}
 \limsup_{n \rightarrow \infty} \Bigl(\frac{1}{2 n \log \log n}\Bigr)^{1/2} \frac{1}{n}\langle u, N_n \rangle & = & 
 -\liminf_{n \rightarrow \infty} \Bigl(\frac{1}{2 n \log \log n}\Bigr)^{1/2} \frac{1}{n}\langle u, N_n \rangle \nonumber \\
 & = & \frac{2}{3\sqrt{3d}}\|u\| \hspace{1cm} \text{a.s.}
\end{eqnarray*}
which immediately leads to
$$
 \limsup_{n \rightarrow \infty}  \Bigl(\frac{1}{2 n \log n \log \log \log n}\Bigr)^{1/2} \frac{1}{n}\langle u, N_n \rangle
 = 0 \hspace{1cm} \text{a.s.}  
$$
On the other hand, we obtain from the law of iterated logarithm for $S_n$ given in Theorem 3.5 of \cite{Bercu19} that 
for any vector $u \in \dR^d$,
\begin{equation*}
\limsup_{n \rightarrow \infty} \Bigl(\frac{1}{2 n \log n \log \log \log n}\Bigr)^{1/2} \langle u, G_n \rangle 
\vspace{-2ex}
\end{equation*}
\begin{eqnarray}
\hspace{1cm}  & = & 
\limsup_{n \rightarrow \infty} \Bigl(\frac{1}{2 n \log n \log \log \log n}\Bigr)^{1/2} \frac{1}{n}\langle u, b_n M_n-N_n \rangle \nonumber \\
& = & 
\limsup_{n \rightarrow \infty} \Bigl(\frac{1}{2 n \log n \log \log \log n}\Bigr)^{1/2} \frac{1}{n}\langle u, b_n M_n\rangle \nonumber \\
& = &
\limsup_{n \rightarrow \infty} \Bigl(\frac{1}{2 n \log n \log \log \log n}\Bigr)^{1/2} \frac{1}{n}\langle u, a_n b_n S_n\rangle \nonumber \\
 &=&
\limsup_{n \rightarrow \infty} \Bigl(\frac{1}{2 n \log n \log \log \log n}\Bigr)^{1/2} \frac{2}{3}\langle u, S_n \rangle \nonumber \\
 &=&
-\liminf_{n \rightarrow \infty} \Bigl(\frac{1}{2 n \log n \log \log \log n}\Bigr)^{1/2} \frac{2}{3}\langle u, S_n \rangle  \nonumber \\
 & = & \frac{2}{3\sqrt{d}}\|u\| \hspace{1cm} \text{a.s.}
 \label{LIL-CM-CR1}
\end{eqnarray}
Hence, we clearly deduce from \eqref{LIL-CM-CR1} that for any vector $u\in\dR^d$,
\begin{equation}
\label{LIL-CM-CR2}
\limsup_{n \rightarrow \infty} \frac{1}{2 n \log n \log \log \log n} \langle u, G_n \rangle^2  = \frac{4}{9d}\|u\|^2 \hspace{1cm} \text{a.s.}
\end{equation}
By taking all rational points on the unit sphere $\dS^{d-1}$ in $\dR^d$, the bound in \eqref{LIL-CM-CR2} holds simultaneously for all of them,
 which implies that
$$
 \limsup_{n \rightarrow \infty} \frac{\| G_n\|^2}{2 n \log n \log \log \log n} \leq   \!\!
 \sup_{u \in \dQ^d \cap \dS^{d-1}} \!\! \limsup_{n \rightarrow \infty} \frac{\langle u, G_n \rangle^2}{2 n \log n \log \log \log n}  =
 \frac{4}{9d} \hspace{1cm} \text{a.s.}
$$
In addition, for any single $u \in \dS^{d-1}$, we also obtain the reverse inequality
$$
 \limsup_{n \rightarrow \infty} \frac{\| G_n\|^2}{2 n \log n \log \log \log n} \geq   
 \limsup_{n \rightarrow \infty} \frac{\langle u, G_n \rangle^2}{2 n \log n \log \log \log n}  =\frac{4}{9d} \hspace{1cm} \text{a.s.}
$$
It immediately leads to \eqref{LIL-CM-CR3} which achieves the proof of Theorem \ref{T-LFQLIL-CM-CR}.
\demend

\subsection{The superdiffusive regime.}
\ \vspace{2ex}\\
\noindent{\bf Proof of Theorem \ref{T-ASCV-CM-SR}.}
It follows from Theorem 3.7 in \cite{Bercu19} that 
\begin{equation}
\label{ASCVG-SN-SR}
\lim_{n \to \infty} \frac{1}{n^a}S_n = L \hspace{1cm} \text{a.s.}
\end{equation}
where the limiting value $L$ is a non-degenerate random vector of $\dR^d$. Hence,
\eqref{ASCVG-SN-SR} together with the Toeplitz lemma imply
\eqref{ASCV-CM-SR1} where the limiting value
$$
G=\frac{1}{a+1}L.
$$
Moreover, we have from \eqref{EMCM} that
\begin{eqnarray*}
 \dE\Big[ \Big\| \frac{1}{n^a}G_n - G \Big\|^2 \Big] & = & \dE\Big[\Big\|\frac{1}{n^{a+1}}(b_nM_n -N_n) - G\Big\|^2\Big], \\
 & \leq & 2\dE\Big[\Big\|\frac{a_n b_n}{n^{a+1}}S_n - G\Big\|^2\Big] + 2\dE\Bigl[\Bigl\| \frac{1}{n^{a+1}}N_n\Bigr\|^2\Bigr].
 \end{eqnarray*}
On the one hand, we already saw from \eqref{CVGANBN} that
$$\lim_{n\to\infty}\frac{a_nb_n}{n}=\frac{1}{a+1}.$$ 
Consequently, we deduce from the mean square convergence (3.12) in \cite{Bercu19} that 
\begin{equation}
\label{MSCVG-CM1}
\lim_{n\to\infty} \dE\Big[\Big\|\frac{a_n b_n}{n^{a+1}}S_n - G\Big\|^2\Big]=0.
\end{equation}
On the other hand, $\dE[\|N_n\|^2] =  \dE[\text{Tr} \langle N_n \rangle)] \leq \tau_n$ where $\tau_n$ is given by \eqref{DEFTAUN}.
Since $\tau_n$ is equivalent to $n^3/3(a+1)^2$ and $a>1/2$, it is not hard to see that
\begin{equation}
\label{MSCVG-CM2}
\lim_{n\to\infty} \dE\Bigl[\Bigl\| \frac{1}{n^{a+1}}N_n\Bigr\|^2\Bigr]= 0.
\end{equation}
Finally, we obtain \eqref{ASCV-CM-SR2} from \eqref{MSCVG-CM1} and \eqref{MSCVG-CM2}, completing the proof of Theorem \ref{T-ASCV-CM-SR}.
\demend

\section*{Appendix C. Proofs of the asymptotic normality results}
\renewcommand{\thesection}{\Alph{section}}
\renewcommand{\theequation}{\thesection.\arabic{equation}}
\setcounter{section}{3}
\setcounter{subsection}{0}
\setcounter{equation}{0}
\setcounter{thm}{0}

\subsection{The diffusive regime.}
\ \vspace{2ex}\\
\noindent{\bf Proof of Theorem \ref{T-CLT-CM-DR}.}
On the one hand, we already saw from \eqref{DECOMPMARTGN} that
\begin{equation*}
\frac{1}{\sqrt{n}}G_n= v^T V_n \cM_n
\hspace{1cm}\text{where}\hspace{1cm}
v=\begin{pmatrix}1 \\ -1\end{pmatrix}\otimes I_d.
\end{equation*}
On the other hand, we deduce from \eqref{CVGIPMN} and \eqref{LINDEBERG-DR4} 
that the two conditions $(\textnormal{H.1})$ and $(\textnormal{H.2})$ of Theorem \ref{T-CLT} are satisfied.
Consequently, we obtain that
\begin{equation*}
\frac{1}{\sqrt{n}}G_n \liml \cN \bigl(0, v^TVv  \bigr)
\end{equation*} 
where the matrix $V$ is given by \eqref{DEFMATRIXV}. It clearly leads to \eqref{CLT-CM-DR1} as
$$
v^TVv= \frac{2}{3(1-2a)(2-a)d}I_d.
$$
\demend
\vspace{-2ex}
\subsection{The critical regime.}
\ \vspace{2ex}\\
\noindent{\bf Proof of Theorem \ref{T-CLT-CM-CR}.}
The proof follows exactly the same lines as that of Theorem \ref{T-CLT-CM-DR} replacing $V_n$ by $W_n$. The details are left to the reader.
\demend



\vspace{-2ex}
\bibliographystyle{abbrv}
\bibliography{Biblio-CDM}
\end{document}